\documentclass[a4paper,11pt, reqno]{amsart} 
\usepackage{amssymb,amsthm,amsmath}

\usepackage{amsfonts}
\usepackage{amscd}
\usepackage{amssymb}
\usepackage{enumerate}
\allowdisplaybreaks
\usepackage{mathabx}

\usepackage{amsbsy}
\usepackage{graphicx}
\usepackage{amsthm}
\usepackage{amsmath}
\usepackage{amsxtra}
\usepackage{mathrsfs}
\usepackage{bbm}

\usepackage{ifthen}

\usepackage[colorlinks,citecolor=red,pagebackref,hypertexnames=false]{hyperref} 

\newtheorem{theorem}{Theorem}[section]
\newtheorem{cor}[theorem]{Corollary}
\newtheorem{lem}[theorem]{Lemma}

\theoremstyle{definition}

\theoremstyle{remark}
\newtheorem{rem}[theorem]{Remark}
\numberwithin{equation}{section}

\newcommand{\R}{\mathbb R}
\newcommand{\T}{\mathbb T}
\newcommand{\Z}{\mathbb Z}
\newcommand{\Na}{\mathbb N}

\newcommand{\rc}{\mathcal{R}}

\newcommand{\ecs}{\mathcal{E}_S}

\newcommand{\la}{\lambda}
\newcommand{\C}{{\mathbb C}}

\newcommand{\pa}{\partial }

\newcommand{\Hc}{\mathcal{H}}

\newcommand{\D}{\Delta}
\newcommand{\Dk}{\Delta_{\kappa}}
\newcommand{\Hk}{H_{\kappa}}
\newcommand{\dl}{\delta }
\newcommand{\ap}{\alpha}
\newcommand{\bt}{\beta }
\newcommand{\N}{\nabla }
\newcommand{\K}{\kappa }
\newcommand{\gm}{\gamma}
\newcommand{\Gm}{\Gamma}
\newcommand{\ve}{\varepsilon}
\newcommand{\gk}{\gamma_\kappa}

\newcommand{\s}{\mathfrak{S}}
\newcommand{\Nk}{\nabla_\kappa}

\renewcommand{\Re}{\operatorname{Re}}

\title[Strichartz inequality for orthonormal functions]
{Strichartz inequality for orthonormal functions associated with Dunkl Laplacian and Hermite-Schr\"{o}dinger operators}

\author[    P Jitendra K. Senapati and Pradeep B]{P Jitendra Kumar Senapati and Pradeep Boggarapu}
\address[P Jitendra K. Senapati.]{Department of Mathematics\\
   BITS Pilani K K Birla Goa Campus\\
    Zuarinagar, South Goa\\
403 726, Goa, India}
\email{jitusnpt@gmail.com}

\address[Pradeep B.]{Department of Mathematics\\
   BITS Pilani K K Birla Goa Campus\\
    Zuarinagar, South Goa\\
403 726, Goa, India}
\email{pradeepb@goa.bits-pilani.ac.in}


\begin{document}

\maketitle
\begin{abstract}
Strichartz inequality for the solutions of free Schr\"odinger equation associated with Dunkl Hermite operator $H_\kappa$ is generalized to any system of orthonormal functions with initial data. A relation between the kernels  of Schr\"odinger propagators ($e^{-it H_\kappa}$ and $e^{it\Delta_\kappa}$) associated with the   Dunkl Hermite and Dunkl Laplacian  operators is established using which corresponding Schtrichartz inequality for orthonormal functions associated with Dunkl Laplacian is obtained.
\end{abstract}

\section{Introduction and main results}
In quantum mechanics, a system of $N$ independent fermions in $\R^d$ can be  described by a collection of $N$ orthonormal functions $f_1,\ldots, f_N$ in $L^2(\R^d).$ For this particular reason, functional inequalities involving a large number of orthonormal functions are very useful in the mathematical analysis of large number of quantum systems. The first such type of generalization goes back to the famous  work by Lieb and Thirring  \cite{LT, LT2}, known as  Lieb-Thirring inequality 
\begin{equation}\label{lti} 
\int_{\R^d}\big(\sum^N_{j=1}|\N f_j(x)|^2\big) dx\geq C \int_{\R^d}\big(\sum^N_{j=1}|f_j(x)|^2\big)^{1+\frac{2}{d}}dx,
\end{equation}  
 which generalizes the  following  Gagliardo-Nirenberg-Sobolev inequality 
\begin{equation}\label{gnsi}
\int_{\R^d}|\N f(x)|^2 dx \geq C'\int_{\R^d}|f(x)|^{2+\frac{4}{d}} dx
\end{equation}
for an $L^2$-normalized function $f$. The inequlaity (\ref{gnsi}) is a  fundamental tool to  understand the stability of mater \cite{HL, HL2, HL3}. 

Let us consider the free Schr\"odinger equation 
\begin{equation}\label{Sch}
\begin{cases}i{u_t}(t, x)&=-\D u(t, x), \hspace{0.5cm}x\in\R^d, t\in\R\\
u(0, x)&=f(x), 
\end{cases}
\end{equation}
where $\D=\sum_{j=1}^d \frac{\pa^2}{\pa x_j^2}$,  Laplacian on $\R^d$. For $f\in L^2(\R^d)$, $e^{it\D}f$ is the unique solution to the initial value problem (\ref{Sch}). The associated Strichartz inequality (see \cite{RS,KY,LS,KT,TT}) reads  as
\begin{equation}\label{Ste1}
\int_{\R}\left(\int_{\R^d}|\left(e^{it\D}f\right)|^{2q}dx\right)^\frac{p}{q}dt\leqslant C\left(\int_{R^d}|f(x)|^2dx\right)^p,
\end{equation}
where $p, q\geqslant1$ satisfy $(p, q, d)\neq(1, \infty, 2)$ and $$\frac{2}{p}+\frac{d}{q}=d.$$
Rupert L. Frank et. al in \cite{FLLS} and \cite{FS} generalized the above inequality for a system of  orthonormal functions and the authors proved the following:

\begin{theorem}\label{SL}
Let $d\geqslant 1$ and $ p, q\geqslant 1$ be such that 
$$\frac{2}{p}+\frac{d}{q}=d,~~~1\leqslant q < \frac{d+1}{d-1}. $$Then for any infinite or finite orthonormal system $(f_j)$ in $L^2(\R^d)$ and for any sequence $(n_j)$ in $\C$, we have 
\begin{equation}\label{Sil}
\int_{\R}\left(\int_{\R^d}\left|\sum_{j} n_j \big | e^{it\D} f_j(x)\big|^2 \right|^q dx\right)^{p/q}dt \leqslant C_{d, q}^p \Big( \sum_j |n_j|^{\frac{2q}{q+1}}\Big)^{\frac{p(q+1)}{2q}}, \end{equation}
 where $C_{d, q}^p$ is a universal constant which only depends on $d$ and $q$.
\end{theorem}
 In \cite{FLLS}, the authors proved that if $q\geqslant \frac{d+1}{d-1}$, then the inequality \eqref{Sil} doesn't hold and the exponent on the right hand side of the inequality \eqref{Sil} is optimal i.e., the exponent $\frac{2q}{q+1}$ can not be increased.

Our aim of this article is to prove the inequality analogues to the inequality \eqref{Sil} when we replace Laplacian by Dunkl Laplacian (for the definition, see Section \ref{Pr}). Motivation to prove inequality of the form \eqref{Sil} in the framework of  Dunkl Laplacian  is due to the recent results of \cite{MJ, MJ2, FLLS, FS, NSS}.  
Let $e^{it\D_k}f$ be the unique solution of the free Schr\"odinger equation associated with Dunkl Laplacian:
\begin{equation}\label{Schd}
\begin{cases}i{u_t}(t, x)&=-\Dk u(t, x), \hspace{0.5cm}x\in\R^d, t\in\R\\
u(0, x)&=f(x). 
\end{cases}
\end{equation}
Mejjaoli studied the above equation first in the literature for the case of Dunkl Laplacian in \cite{M1} and he proved the following Strichartz inequality:
\begin{equation}\label{Sted}
\int_{\R}\left(\int_{\R^d}|e^{it\Dk}f(x)|^{2q}dx\right)^\frac{p}{q}dt\leqslant C\left(\int_{R^d}|f(x)|^2dx\right)^p,
\end{equation}
where $d\geqslant 1$ and $p, q\geq 1$ are such that 
$$\frac{2}{p}+\frac{d+2\gk}{q}=d+2\gk.$$
As applications of study of the Strichartz estimates \eqref{Sted} for the Dunkl Schr\"{o}dinger equation \eqref{Schd}, Mejjaoli had studied the global well-posedness and scattering for a class of non linear Dunkl Schr\"{o}dinger equations in \cite{M2}. For all unexplained notations, we refer to Section \ref{Pr}.  In the same spirit of \cite{FLLS} and \cite{FS}, it is natural to ask the inequality analogues to \eqref{Sil} for the case Dunkl Laplacian. Our main result lies on this direction:

\begin{theorem}[Strichartz inequality for orthonormal functions for Dunkl Laplacian] \label{TSOF}
Let $d\geqslant 1$ and $p,q\geqslant1$ satisfy $$1\leqslant q < {\frac{d+2\gk+1}{d+2\gk-1}}~ \mbox{and}~ \frac{2}{p}+\frac{d+2\gk}{q}=d+2\gk.$$
Then for any (possibly infinite) system $(f_j)$ of orthonormal functions in $L_{\kappa}^2 (\R^d)$ and any coefficients $(n_j)\subset{\C},$ we have
\begin{equation}\label{Ste2} \int_\R\left(\int_{\R^d}\left|\sum_{j}n_j|e^{it\Dk}f_j(x)|^2\right|^q dx\right)^{\frac{p}{q}}dt\leqslant{{C^{p}_{d,q}}\left(\sum_{j}|n_j|^{\frac{2q}{q+1}}\right)}^\frac{p(q+1)}{2q},
\end{equation}
where $C_{d,q}$ is a universal constant which only depends on $d$ and $q.$
\end{theorem}

\begin{rem}
For $q=1$ and $p=\infty,$ since $e^{it\D_\kappa}$ is a unitary operator on $L^2_\kappa(\R^d),$  using triangle inequality it can be seen that $$\sup_{t\in\R}\left(\int_{\R^d}\left|\sum_j n_j|e^{it\D_\kappa}f_j(x)|^2\right|h^2_\kappa(x)dx\right)\leqslant\sum_j|n_j|.$$
\end{rem}

Further extension of inequality \eqref{Ste1} has been made for the Schr\"odinger equation of the form $iu_t(t, x)-P u(t, x) - V(x)u(t, x) = 0,$ for a suitable potential $V$ by several authors, for example, see \cite{GV, JSS, KPV} when $P=-\D$ and \cite{M3, M4} when $P=-\Dk$. In particular, when $V(x)=|x|^2$,  the Strichartz inequalities have been studied in the literature see \cite{KT, NR, MJ}. It is natural to ask whether such inequalities are true when  $\D$ is replaced by $\Dk$. In this case the initial value problem \eqref{Sch} turns out to be an initial value problem for the Schr\"odinger equation associated with the Dunkl Hermite operator $\Hk=-\Dk+|x|^2$:  
\begin{equation}\label{Schdh} 
\begin{cases}i{u_t}(t, x)&=H_\K u(t, x), \hspace{0.5cm}x\in\R^d, t\in\R,\\
u(0, x)&=f(x). 
\end{cases}
\end{equation}
 The author in \cite{M3} proved the Strichartz inequalities for the above problem. If $f\in L^2_\K(\R^d),$ the solution of the initial value problem \eqref{Schdh} is given by $u(t, x) = e^{-itH_\K}f(x).$ The Strichartz inequality in this case reads as:

\begin{theorem}\cite{M3} Let $f \in L^2_\K(\R^d).$ If $p, q \geq 1$ satisfying 
$$\left( \frac{d+2\gk-2}{d+2\gk}\right)<\frac{1}{q}\leq 1 ~~\text{and}~~1\leq \frac{1}{p} \leq 2, $$ or $$ 0\leq \frac{1}{p} <1 ~~\text{and}~~ \frac{2}{p}+\frac{d+2\gk}{q}\geq (d+2\gk).$$
Then 
\begin{equation}\label{Stedh}
\|e^{-itH_\K}f\|_{L^{2p}_t L^{2q}_{\K, x}(\T \times \R^d)} \leq C \|f\|_{L^2_\K(\R^d)},
\end{equation}
where $\T=(-\pi/4, \pi/4), (-\pi/2, \pi/2)$ or $(-\pi, \pi)$.
\end{theorem}
We are interested to extend the inequality \eqref{Stedh} for orthonormal family of functions and  prove the following :
\begin{theorem}\label{TSOFdh}
Let $d\geqslant 1$ and $p,q\geqslant1$ satisfy $$1\leqslant q < {\frac{d+2\gk+1}{d+2\gk-1}}~~ \mbox{and}~~ \frac{2}{p}+\frac{d+2\gk}{q}=d+2\gk.$$
Then for any (possibly infinite) system $(f_j)$ of orthonormal functions in $L_{\kappa}^2 (\R^d)$ and any coefficients $(n_j)\subset{\C},$ we have
\begin{equation}\label{Ste3} \int^\pi_{-\pi}\left(\int_{\R^d}\left|\sum_{j}n_j|\left(e^{-it\Hk}f_j\right)(x)|^2\right|^q dx\right)^{\frac{p}{q}}dt\leqslant{{C^{p}_{d,q}}\left(\sum_{j}|n_j|^{\frac{2q}{q+1}}\right)}^\frac{p(q+1)}{2q},
\end{equation}
where $C_{d,q}$ is a universal constant which only depends on $d$ and $q.$
\end{theorem}
\begin{rem}\label{rem16}
For $q=1$ and $p=\infty,$ since $e^{-it\Hk}$ is a unitary operator on $L^2_\kappa(\R^d),$  using triangle inequality it can be seen that $$\sup_{t\in(-\pi, \pi)}\left(\int_{\R^d}\left|\sum_j n_j|e^{-it\Hk}f_j(x)|^2\right|h^2_\kappa(x)dx\right)\leqslant\sum_j|n_j|.$$
\end{rem}
In \cite{RS, FS}, the authors proved  Fourier restriction theorems and they linked these theorems to space time decay estimates for certain evolution equations. We prove Theorem \ref{TSOF} assuming Theorem \ref{TSOFdh} and the relation between the Schr\"odinger kernels associated with $\Hk$ and $\Dk$, see Lemma \ref{rel}. In proving Theorem \ref{TSOFdh}, we closely follow the ideas of famous work of Frank-Sabin \cite{FS} and Strichartz \cite{RS}. For better understanding of the ideas we need to start with Fourier-Dunkl Hermite transform.

Let $f\in L^1_\K((-\pi, \pi)\times \R^d)$ and $\Na_0$ denotes the set of all non negative integers.  We define the Fourier-Dunkl Hermite transforms of $f$ by $$ \hat{f}(\mu, \nu) :=\int_{-\pi}^{\pi} \int_{\R^d} f(t, x)\phi^\K_\mu(x) e^{i\nu t} h^2_\K(x) dx dt,~~~ ~~(\mu, \nu)\in\Na_0^d\times \Z,$$
where $\phi^\K_\mu$'s are the $d$-dimensional Dunkl Hermite functions (see Section \ref{Pr}). Since the set $\{\phi_\mu^\K: \mu \in \Na_0^d\}$ is a complete orthonormal system for $L^2_\K(\R^d)$, if $f\in  L^2_\K((-\pi, \pi)\times \R^d)$ then $\{\hat{f}(\mu, \nu)\} \in \ell^2(\Na_0^d \times \mathbb{Z})$ and satisfies the Plancheral formula $$\|f\|^2_{ L^2_\K((-\pi, \pi)\times \R^d)}=\sum_{\mu, \nu}|\hat{f}(\mu, \nu)|^2.$$ The inverse Fourier-Dunkl Hermite transform is given by $$f(t, x) = \sum_{\mu, \nu} \hat{f} (\mu, \nu) \phi^\K_\mu(x) e^{-i\nu t}.$$   

Given a discrete surface $S$ in $\Na^d_0\times \Z,$ we define the restriction operator on suitable functions by $\rc_s f :=\{{\hat{f} (\mu, \nu)}\}_{(\mu, \nu)\in S}$ and the operator dual to $\rc_S$ is called the extension operator defined as $$\ecs(\{g(\mu, \nu)\}) := \sum_{(\mu, \nu) \in S} g(\mu, \nu) \phi^\K_\mu(x)e^{-i\nu t}$$ and ask the following question. For which exponents $1\leq p \leq 2,$  Fourier-Dunkl Hermite transforms of a function $f \in L^p_\K((-\pi, \pi)\times \R^d)$ belongs to $\ell^2(S)?$ This question can be reframed to the boundedness of the operator $\ecs$ from $\ell^2(S)$ to $L^{p'}_\K((-\pi, \pi)\times \R^d),$ where $p'$ is the conjugate exponent of $p$ by a duality argument. Since $\ecs$ is bounded from $\ell^2(S)~ \mbox{to}~ L^{p'}_\K((-\pi, \pi)\times \R^d)$ if and only if $T_S := \ecs (\ecs)^*$ is bounded from $L^p_\K((-\pi, \pi)\times \R^d)~ \mbox{to}~ L^{p'}_\K((-\pi, \pi)\times \R^d)$. Theorem of this problem in the literature is known as restriction theorem.

Restriction theorems for the Fourier transforms for certain quadratic surfaces $S$ of $\R^N$ have been proved by Stein \cite{ES} and Strichartz \cite{RS} using the Stein's interpolation theorem \cite{ES2}.  As an application of H\"older's inequality, it can be verified that $T_S$ is $L^p-L^{p'}$ bounded if and only if the operator $W_1T_S W_2$ is bounded from $L^2$ to $L^2$ for any $W_1, W_2\in L^{\frac{2p}{2-p}}$. Recently, Frank and Sabin \cite{FS} considered the similar problem for the case of Fourier transform and they proved that for any $W_1, W_2 \in  L^{\frac{2p}{2-p}}(\R^N)$, the operator $W_1T_S W_2$ is not only bounded operator on $L^2$, also belongs to a Schatten class (for definitions, see Subsection \ref{Sbss}). The estimate that they proved is
$$\|W_1 T_S W_2\|_{\s^\ap(L^2{(\R^N)})}\leq C \|W_1\|_{L^{\frac{2p}{2-p}}(\R^N)}  \|W_2\|_{L^{\frac{2p}{2-p}}(\R^N)},$$ for some $\alpha \geq 1$, where $C>0$ independent of $W_1$ and $W_2$. Using these kind of estimates along with duality argument (Lemma 3 in \cite{FS}), they proved inequality \eqref{Sil}.

We consider a particular discrete surface $S\subset \Na_0^d \times \Z$  with respect to counting measure for which $\ecs=e^{-it\Hk}$ and define an analytic family of operators $(T_z)$ defined on the strip $a \leq \text{Re} ~z \leq b$ in the complex plane such that $T_S=T_c$ for some $c\in (a, b)$. Then we show that the operator $W_1T_S W_2$ belongs to a Schatten class for $W_1, W_2$ in suitable mixed norm spaces and applying the duality argument (see Lemma \ref{Ld}) we obtain the Strichartz inequality \eqref{Ste3} for the system of orthonormal functions for Dunkl Hermite operator $\Hk$. In order to prove Theorem \ref{TSOF}, we will establish a relation between Schr\"odinger kernels associated with Dunkl Laplacian and Dunkl Hermite operator using which along with Theorem \ref{TSOFdh} we obtain the inequality \eqref{Ste2}. 



%
%

\subsection{Dual Strichartz inequality.}
In this subsection, we prove that  inequalities \eqref{Ste2} and \eqref{Ste3} are  equivalent to some dual inequalities.
One can verify that Theorem \ref{TSOFdh} (or Theorem \ref{TSOF}) is equivalent to the boundedness of the operator $\gm \mapsto \rho_{\gamma(t)}$ from $\s^{\frac{2q}{q+1}}(L^2_\K(\R^d))$ to $L^p_t L^q_{\kappa, x}(\T \times\R^d)$, where $\gm(t)=e^{-itP}\gm e^{itP}$ with $P=\Hk$ and $\T=(-\pi, \pi)$ (or $P=-\Dk$ and $\T=\R$) (See Subsection \ref{Sbss}). Dual of this operator is given by 
$$V\in L^{p'}_t L^{q'}_{\K, x}(\T \times \R^d) \mapsto\int_\T e^{itP}V(t, x)e^{-itP}dt\in\s^{2q'}(L^2_\K(\R^d)).$$

Therefore the duality argument shows that Theorem \ref{TSOFdh} and Theorem \ref{TSOF} are equivalent to the following:

\begin{theorem}[Strichartz inequality in Schatten spaces, dual version] \label{TSOFD}
Assume that $p', q', d\geq1$ satisfy $$1+\frac{d+2\gk}{2}\leq q'<\infty~ \mbox{and}~ \frac{2}{p'}+\frac{d+2\gk}{q'}=2.$$
We have
\begin{equation}\label{Sc}
\left\| \int_\T e^{itP} V(t, x) e^{-itP} dt \right\|_{\s^{2q'}(L^2_\K(\R^d))}\leq C_{d, q}\|V\|_{L^{p'}_t L_{\K, x}^{q'}(\T \times\R^d))},
\end{equation}
where $C_{d, q}$ is the same constant as in Theorem \ref{TSOF} when $P=-\Dk, \T = \R$ or as in Theorem \ref{TSOFdh} when $P=\Hk, \T = (-\pi, \pi)$.
\end{theorem}

\begin{rem}

For $q' = \infty,$ and $p' = 1,$ it can be easily seen that 
$$\left\| \int_\T e^{itP} V(t, x) e^{-itP} dt \right\| \leq \|V\|_{L^{1}_t L_{\K, x}^{\infty}(\T \times\R^d)}.$$
\end{rem}

Using \eqref{Sc} when $P=\Hk$, we prove an inhomogeneous Strichartz inequality. Consider the equation given by
\begin{equation}\label{Inh}  
\begin{cases}
i\dot{\gm}(t)&= [-H_\K, \gm(t)] + iR(t)\\
\gm(t_0)&=0,
\end{cases}
\end{equation}
where $R(t)$ is a self adjoint operator on $L^2_\K(\R^d)$ and is bounded for almost every $t.$ The solution of the system \eqref{Inh} can be written as 
\begin{equation}\label{Inh2}
\gm(t) = \int^t_{t_0} e^{i(t-s)H_\K} R(s)  e^{i(s-t)H_\K} ds.
\end{equation}
We obtain the following inhomogeneous Strichartz inequality.
\begin{cor}[Inhomogeneous Strichartz inequality]\label{ISI} 
Assume that $p, q, d \geq 1$ satisfy 
$$1 \leq q < \frac{d + 2\gm_\K+1}{d+2\gm_\K-1}~ \mbox{and}~ \frac{2}{p}+\frac{d+2\gm_\K}{q}=d+2\gm_\K$$
and let $\gm(t)$ be given by \eqref{Inh2}. Then 
$$\|\rho_{\gm(t)}\|_{L^p_tL^q_{\K, x}((-\pi, \pi)\times\R^d)} \leq C \left\|\int^\pi_{-\pi} e^{isH_\K} |R(s)| e^{-isH_\K} ds\right\|_{\s^\frac{2q}{q+1}},$$
for a constant $C$ which is independent of $t_0.$
\end{cor}
\begin{proof} For a trace-class operator $\gamma$ and any bounded function $V$ of compact support, 
$$Tr(V(x)\gamma) = \int_{\R^d} V(x) \rho_{\gm}(x) h_\K^2(x)dx,$$ where $V(x)$ is identified with the corresponding multiplication operator on $L^2_\K(\R^d)$ and $\rho_\gm (x)$ is the density of $\gm$. And for a time-dependent potential $V(t, x) \in L^\infty_c((-\pi, \pi) \times \R^d),$ we have
 \begin{align*}
&\left|\int^\pi_{t_0}\int_{\R^d} V(t, x) \rho_{\gm(t)}(x)  h_\K^2(x) dx dt\right|\\ &= \left|\int^\pi_{t_0}\int^t_{t_0} Tr(e^{itH_\K} V(t, x) e^{-itH_\K} e^{isH_\K} R(s) e^{-isH_\K}) ds dt\right|\\
&\leq \int^\pi_{t_0}\int^t_{t_0} Tr(e^{itH_\K} |V(t, x)| e^{-itH_\K} e^{isH_\K} |R(s)|\ e^{-isH_\K}) ds dt\\
&\leq Tr \left(\left( \int^\pi_{t_0} e^{itH_\K} |V(t, x)| e^{-itH_\K}\right) \left(\int^\pi_{t_0}e^{isH_\K} |R(s)|\ e^{-isH_\K}\right)\right).
\end{align*}
In the first inequality we have used the fact that $|Tr(AB)| \leq Tr(|A||B|)$ for self-adjoint operators $A$ and $B.$ Then applying the H\"older's inequality for traces and \eqref{Sc} for the term involving $V(t, x).$ This completes the proof.      
\end{proof}

%
%
%


We end this section with an overview of our paper. In Section 2, we recall the basic ingredients in Dunkl setting including Dunkl Hermite functions and their generating function identity and  Schatten class of operators. We also state duality statements. We will establish necessary theorems and proposition to give the proofs of our main results Theorem \ref{TSOF}, \ref{TSOFdh}, and \ref{TSOFD} in Section 3. As an application of inhomogeneous Strichartz inequality, we prove the global well-posedness for the Dunkl Hermite-Hartree equation in Schatten spaces in Section 4. In Section 5, we will give an alternate proof of Theorem \ref{TSOF} for $1<q \leq 1+\frac{2}{d+ 2\gk}$.


\section{Preliminaries}\label{Pr}
~~~In this section we will define Dunkl operators, some other related operators and function spaces which will be used in this article. Dunkl operators were introduced by Charles Dunkl (1989) to built a framework for a theory of special functions and integral transforms in several variables related to reflection groups. Such operators are relevant in physics, namely for the analysis of quantum many body systems of Calogero-Moser-Sutherland type (see \cite{DV,LV}). From the mathematical analysis point of view, the importance of Dunkl operators lies on the fact that they generalize the theory of symmetric spaces of Euclidean type. There are many developments in harmonic analysis of the operators which are defined in terms of Dunkl operators in different directions in recent years.  There is a vast literature related to Dunkl operators, see for instance \cite{BRS,CFD,CFD2,DX,H,MR}.

\subsection{The general Dunkl setting}
~~~The basic ingredients in the theory of Dunkl operators are the root systems and finite reflection groups associated to them. For $\nu\in\R^d\setminus\{0\}$, we denote by $\sigma _\nu$ the reflection in the hyperplane perpendicular to $\nu,$ i.e.,
$$\sigma _\nu(x)=x-2\frac{\langle\nu, x\rangle}{|\nu|^2}\nu.$$ 

 Let $O(d)$ be the group of orthogonal matrices acting on $\R^d.$ Given a root system $R$ associate a finite subgroup $G\subset O(d)$ is the reflection group which is generated by the reflections $\{\sigma_{\nu}: \nu\in R\}.$

A function $\kappa : R\rightarrow \C$ is said to be a multiplicity function on $R,$ if it is invariant under the natural action of $G$ on $R,$ i.e. $\kappa(g\nu)=\kappa(\nu) $ for all $\nu\in R$ and $g\in G.$ 

Every root system can be written as a disjoint union $R = R_+ \cup (-R_+)$, where $R_+$ and $-R_+$ are separated by a hyperplane through the origin. Such $R_+$ is called the set of all positive roots in $R$. Of course its choice is not unique.

The weight function associated to the root system $R$ and the multiplicity function $\kappa$ is defined by 
$$h^2_{k}(x):=\prod_{\nu\in R_+}|\langle x, \nu\rangle|^{2\kappa(\nu)}.$$

Note that $h^2_{\kappa}(x)$ is $G$-invariant and homogeneous of degree $2\gamma_\kappa$ where, by definition,
$$\gamma_\kappa:=\sum_{\nu\in R_{+}}\kappa(\nu).$$

We will assume through out the article $2\gk$ is a non negative integer. Let $L^p_{\kappa}(\R^d)$, $1\leq p\leq \infty,$  stands for the space of $L^p-$functions with respect to the measure $h_\kappa^2(x) dx$, and $L^p_t L^q_{\kappa, x}(\T \times \R^d), 1\leq p, q \leq \infty$ stands for the space of all measurable functions $h(t, x)$ on $\T  \times \R^d$   for which
$$\|h\|_{L^p_t L^q_{\kappa, x}(\T \times \R^d)}:=\|\|h(t, \cdot)\|_{L^q_{\kappa}(\R^d)}\|_{L^p(\T)}<\infty.$$
We may consider $\T=(-\pi, \pi), \R$ or any interval in $\R$ with Lebesgue measure. $L^p_\K(\T\times \R^d):=L^p_t L^p_{\kappa, x}(\T \times \R^d)$.

Now we define the difference-differential operators which were introduced and studied by C. F. Dunkl (for $\kappa\geq 0$), see \cite{CFD,CFD2}. These operators are also called Dunkl operators and are analogues (generalization) of directional derivatives. We fix a root system $R$ with a positive subsystem $R_+$ and the associated reflection group $G.$ We also fix a non-negative multiplicity function $\kappa$ defined on $R.$

For $\xi\in\R^d,$ the Dunkl operator $T_{\xi}:=T_{\xi}(\kappa)$ is defined by 
$$T_{\xi}f(x)=\partial_{\xi}f(x)+\sum_{\nu\in R_{+}}\kappa(\nu)\langle\nu, \xi\rangle\frac{f(x)-f(\sigma_{\nu}x)}{\langle\nu, x\rangle},$$
for smooth functions $f$ on $\R^d.$ Here $\partial_{\xi}$ denotes the directional derivative along $\xi.$ For the standard coordinate vectors $\xi=e_{j}$ of $\R^d$ we use the abbreviation $T_{j}=T_{e_j}.$

Let $\mathcal{P}$ be the space of all polynomials with complex coefficients in $d$-variables and $\mathcal{P}_m$ the subspace of homogeneous polynomials of degree $m$. The Dunkl-operators $T_\xi$ and directional derivatives $\partial_\xi$ are closely related and intertwined by an isomorphism on $\mathcal{P}.$ indeed, if the multiplicity function $\kappa$ is non-negative then by Theorem 2.3 and Proposition 2.3 in R\"osler \cite{MR}, there exist a unique linear isomorphism (intertwining operator) $V_\kappa$ of $\mathcal{P}$ such that $V_\kappa(\mathcal{P}_m)=\mathcal{P}_m, V_\kappa|_{\mathcal{P}_0}=id$ and $T_{\xi} V_\kappa=V_\kappa\partial_\xi$ for all $\xi\in\R^d.$ It can be checked that $V_\kappa\circ g=g\circ V_\kappa$ for all $g\in G.$ 

For $y\in\C^d,$ define 
$$E_\kappa(x, y):=V_\kappa\left(e^{\langle \cdot, y\rangle}\right)(x),~ x\in\R^d.$$
The function $E_\kappa$ is called the Dunkl-kernel, or $\kappa$-exponential kernel, associated with $G$ and $\kappa,$ see \cite{CFD4}. 
The Dunkl-Laplacian is the second order operator defined by 
$$\Delta_{\kappa}=\sum^d_{j=1}T^2_j,$$ 
which can be explicitly calculated, see Theorem 4.4.9 in Dunkl-Xu \cite{DX}. 

It can be seen that $\Delta_\kappa=\sum^d_{j=1}T^2_{\xi_j}$ for any orthonormal basis $\{\xi_1, \xi_2,\ldots, \xi_d\}$ of $\R^d$, see \cite{CFD2}.
We also recall the definition of heat kernel, see \cite{MR2}. We are considering the initial value problem for heat equation associated to the Dunkl Laplacian,   
$$u_t(t, x)=\Dk u(t, x), ~~~~ u(0, x) = f(x), ~~~ t>0,~~ x\in \R^d.$$
For $f \in L^p_\K(\R^d)$ the solution of this equation is given by
$$u(t, x) =e^{t\Dk}f(x)= \int_{\R^d} \Gamma_\K(t, x, y) f(y) h^2_\K(y) dy, $$
where the heat kernel $\Gamma_\K$ associated to the Dunkl Laplacian is explicitly given by  

\begin{equation}\label{l4}
\Gm_\K(t, x, y) =\frac{M_\K}{ (2t)^{\gk+d/2}}e^{-\frac{|x|^2+|y|^2}{4t}}E_\K\Big(\frac{x}{2t}, y\Big),
\end{equation}
where $\displaystyle M_\K=\left(\int_{\R^d} e^{-\frac{|x|^2}{2}} h_\K^2(x) dx \right)^{-1}$.
The kernel of the operator $e^{it\Dk}$ is given by
\begin{equation}\label{lsk} L_{it}(x, y)=\frac{M_\K}{ (2 i t)^{\gk+d/2}}e^{\frac{i(|x|^2+|y|^2)}{4t}}E_\K\Big(\frac{x}{2it}, y\Big).\end{equation}

\subsection{Dunkl harmonic oscillator.}
The Dunkl harmonic oscillator (which we also call the Dunkl Hermite operator) is defined by $$H_\K = -\Dk + |x|^2.$$ 
Let $\Na_0$ be the set of all non-negative integers and for the multi index $\mu = (\mu_1, . . ., \mu_d) \in \Na^d_0$ R\"osler introduced eigenfunctions $\phi^\K_\mu$ of Dunkl harmonic oscillator $H_\K$ in \cite{MR2} with eigenvalue $(2|\mu| + d + 2\gk)$, where $|\mu| = \mu_1 + \cdot \cdot \cdot + \mu_d,$ i. e.,
$$H_\K \phi^\K_\mu = (2|\mu| + d + 2\gk) \phi^\K_\mu,$$
and they form a complete orthonormal system in $L^2_\K(\R^d).$  These functions $\phi^\K_\mu$ are called the Dunkel Hermite functions (also called generalized Hermite functions) and they satisfy the following generating function identity:

\begin{equation}\label{Mhef}
\sum_{\mu\in\Na^d_0} \phi^\K_\mu(x) \phi^\K_\mu(y) w^{|\mu|}= 2^{\frac{d}{2}+\gk} M_\K (1-w^2)^{-(\frac{d}{2}+\gk)} e^{-\frac{1}{2}\left(\frac{1+w^2}{1-w^2}\right)(|x|^2+|y|^2)} E_\K\left(\frac{2w x}{1-w^2}, y\right).
\end{equation}

For any $\phi \in L^2_\K(\R^d)$, $\displaystyle \phi=\sum_{\mu\in \Na_0^d} \hat{\phi}(\mu) \phi_\mu^\K$, where 
$$\hat{\phi}(\mu)=\int_{\R^d}\phi(x) \phi_\mu^\K (x) h_\K^2(x) dx.$$  
 For the $Re(z)>0$, using spectral theory
$$e^{-z\Hk}f(x)=\int_{\R^d} K_z(x, y)f(y) h_\K^2(y) dy,$$ 
for $f\in L^2_\K(\R^d)$, where  the kernel 
\begin{align*}
K_z(x, y) &= \sum_{\mu\in \Na^d_0} e^{-z(2|\mu|+d+2\gk)} \phi^\K_\mu(x) \phi^\K_\mu(y)\\
&=e^{-(d+2\gk)z}\sum_{\mu\in \Na^d_0} e^{-2|\mu|z} \phi^\K_\mu(x) \phi^\K_\mu(y).
\end{align*} 
In the view of Mehler's formula \eqref{Mhef} for $w=e^{-2z}$, we will get 

\begin{align*}
K_z(x, y) &=\frac{M_\K}{(\sinh 2z)^{\frac{d}{2}+\gk}} e^{-\frac{1}{2} \coth 2z(|x|^2+|y|^2)}E_\K\left(\frac{x}{\sinh 2z}, y\right).
\end{align*}

Let $z = r+it$ with $r\rightarrow 0^+$, we will get the kernel of the operator $e^{-it\Hk}$:
\begin{equation}\label{hsk}
K_{it}(x,y) = \frac{M_\K}{( i \sin 2t)^{\frac{d}{2}+\gk}} e^{\frac{i}{2} \cot 2t(|x|^2+|y|^2)}E_\K\left(\frac{x}{i \sin 2t}, y\right),
\end{equation}
for all $t\in \R\setminus (\frac{\pi}{2})\mathbb{Z}$. It can be easily verified that the following:
\begin{equation}\label{pk1} 
K_{-it}(x, y)=\overline{K_{it}(x, y)}~~~\quad~~\text{and}~~~\quad~~~K_{i(t+\frac{\pi}{2})}(x, y)=(-1)^{\frac{d}{2}+\gk}K_{it}(-x, y),
\end{equation}
for all $t\in \R\setminus (\frac{\pi}{2})\mathbb{Z}$.


\subsection{Schatten Spaces.}\label{Sbss}
Let $\Hc$ be complex separable Hilbert space and $\gamma$ be a compact operator on $\Hc$.  We say that $\gamma\in \s^p(\Hc),$ Schatten space, for $ 1\leqslant p<\infty,$ if $Tr|\gamma|^p<\infty,$ where $|\gamma|=\sqrt{\gamma^*\gamma}$ and $Tr(\gamma)$ is trace of $\gamma$. For $\gamma\in\s^p(\Hc),$ the Schatten $p$-norm of $\gamma$ in $\s^p(\Hc)$ is defined by 
$$\|\gamma\|_{\s^p(\Hc)}=(Tr|\gamma|^p)^\frac{1}{p}.$$
It can be verified that  $\|\gamma\|_{\s^p(\Hc)}=\left(\sum_{j}|\lambda_{j}|^p\right)^\frac{1}{p},$ where $\lambda_1\geqslant\lambda_2\geqslant \dots \geqslant\lambda_n\geqslant \dots \geqslant 0$ are the singular values of $\gamma,$ that is the eigenvalues of $|\gamma|=\sqrt{\gamma^*\gamma}.$

When $p=2,$  the Schatten $p$-norm coincides with the Hilbert-Schmidt norm.  Also when $p=\infty,$ we define $\|\gamma\|_{\s^\infty(\Hc)}$ to be the operator norm of $\gamma$ on $\Hc.$ For $\gm \in \s^p(L^2_\K(\R^d)), \rho_\gm$ denotes the density of $\gm$ and it is defined for any finite rank $\gm$ by duality,
$$\int_{\R^d} \rho_\gm (x) V(x) h_\K^2(x) dx := Tr(\gm V)$$ and extended to all $\gm \in \s^p(L^2_\K(\R^d))$ using the density of finite rank operators in $\s^p(L^2_\K(\R^d)$. We refer \cite{S} for more details of Schatten spaces.

 We shall write $C$ to denote positive constants independent of significant quantities the meaning of which can change from one occurrence to another. The following lemma plays an important role in proving the main results.

\begin{lem}[Duality principle]\label{Ld}
Let $p, q \geq 1~ \mbox{and} ~ \ap \geq 1.$  Let $A$ be a bounded linear operator from $L^2_\K(\R^d)$ to $L^{2p}_t L^{2q}_{\K, x}(\T \times \R^d)$. Then the following statements are equivalent.

$(1)$ There is a constant $C > 0$ such that 
\begin{equation}\label{DI}
\|WAA^*\overline{W}\|_{\s^\ap\left(L^2_{\K}(\T\times \R^d)\right)} \leq C \|W\|^2_{L^{2p'}_t L^{2q'}_{\K, x}(\T \times \R^d)},
\end{equation}
for all $W\in L^{2p'}_t L^{2q'}_{\K, x}(\T \times \R^d),$ where the function $W$ is interpreted as an operator which acts by multiplication.\\

$(2)$ For any orthonormal system $(f_j)_{j\in J}$ in $\Hc$ and any sequence $(n_j)_{j\in J}\subset \C,$ 

\begin{equation}\label{DI2}
\left\|\sum_{j\in J}n_j\left|Af_j\right|^2\right\|_{L^p_t L^q_{\K, x}(\T \times \R^d)}\leq C' \left(\sum_{j\in J}|n_j|^{\ap'}\right)^{1/\ap'},
\end{equation}
where $C'$ is a constant. 
\end{lem}
For the proof of above lemma we refer Lemma 3 in \cite{FS} with appropriate modification. We end this section with following two important remarks.

\begin{rem}
Let $S$ be the discrete surface $S = \{(\mu, \nu) \in \Na^d_0 \times \Z : \nu = 2|\mu|+d+2\gk \}$ with respect to the counting measure. Then for all $\{g(\mu, \nu)\} \in L^1(S)$ and for all $(t, x) \in (-\pi, \pi)\times\R^d,$  the extension operator can be written as  
\begin{equation}\label{Eoe}
\ecs(\{g(\mu, \nu)\}) (t, x) = \sum_{\mu, \nu \in S} g (\mu, \nu) \phi^\K_\mu(x) e^{-it\nu}.
\end{equation} Choosing
\begin{equation} \label{cg}
g(\mu, \nu) = \begin{cases}\hat{u}(\mu) & \text{if}~ \nu=2|\mu|+d+2\gk,\\
0 &\text{otherwise}, 
\end{cases} 
\end{equation}
for any suitable $u : \R^d \rightarrow \C$, where $\hat{u}(\mu)=\int_{\R^d} u(x) \phi_\mu^\K (x) h_\K^2(x) dx$. We will get 
\begin{align*}
\ecs(\{g(\mu, \nu)\}) (t, x)  &= \sum_{\mu, \nu \in S} g (\mu, \nu) \phi^\K_\mu(x) e^{-it\nu}\\
& = \int_{\R^d} \left(\sum_\mu \phi^\K_\mu (x) \phi^\K_\mu (y) e^{-it(2|\mu|+d+2\gk)} \right) u(y) h^2_\K(y) dy\\
& =e^{-itH_\K} u(x).
\end{align*}
Thus we have for $S = \{(\mu, \nu) \in \Na^d_0 \times \Z : \nu = 2|\mu|+d+2\gk \}$,
\begin{equation}\label{210}
\ecs(\{g(\mu, \nu)\}) (t, x) = e^{-itH_\K} u(x),
\end{equation}
whenever $g$ and $u$ are related by \eqref{cg}. 
\end{rem}

\begin{rem}
For any compact operator $\gamma$ on $L^2_\kappa(\R^d)$, we define
$$\gamma(t) :=e^{-itP}\gamma e^{itP}, ~~t\in\R,$$where $P=-\Dk$ or $\Hk$. If the operator $\gamma$ is of the form $\gamma :=\sum_{j}n_j|f_j\rangle\langle f_j|$ associated with a given orthonormal system $(f_j)_j,$ where $|f \rangle \langle g |$ is Dirac's notation for the rank one operator $\phi\mapsto\langle\phi, g\rangle f$,  for such $\gamma$ one may check that 
$$\rho_{\gamma(t)}(x) :=\sum_jn_j|e^{-itP}f_j(x)|^2.$$
Using this relation, the inequality \eqref{Ste3} in Theorem \ref{TSOFdh} (or the inequality \eqref{Ste2} in Theorem \ref{TSOF})   can also be written as
\begin{equation}\label{op1} 
\|\rho_{\gamma(t)}(x)\|_{L^p_tL^q_{\K, x}(\T \times \R^d)}\leqslant C_{d, q}\|\gamma\|_{\s^{\frac{2q}{q+1}}(L^2_\K(\R^d))},
\end{equation}
where $\|\gamma\|_{\s^{\frac{2q}{q+1}}(L^2_\K(\R^d))} = {\left(\sum_{j}|n_j|^{\frac{2q}{q+1}}\right)}^\frac{q+1}{2q}, P=\Hk$ and $\T =(-\pi, \pi)$ (or $P=-\Dk$ and $\T=\R$). The estimate \eqref{op1} is equivalent to say that $\gm \mapsto \rho_{\gm(t)}$ is a bounded operator from $\s^{\frac{2q}{q+1}}(L^2_\K(\R^d))$ to $L^p_tL^q_{\K, x}(\T \times \R^d)$.

 Similarly in view Lemma \ref{Ld}, the estimates \eqref{DI} and  \eqref{DI2} for $A=e^{-itP}$ are equivalent to the following bound: for any $\gamma\in\s^{\ap'}(L^2_{\K}(\R^d))$ we have 
$$\|\rho_{e^{-itP}\gamma e^{itp}}\|_{L^p_t L^q_{\K, x}(\T \times \R^d)} \leq C \|\gm\|_{\gamma\in\s^{\ap'}\left(L^2_{\K}(\R^d)\right)}.$$    
\end{rem}

\section{Proofs}

\subsection{Proof of Theorem \ref{TSOF}}\label{pt1}
In this subsection we will give the proof of Theorem \ref{TSOF} by assuming Theorem \ref{TSOFdh} and we will prove Theorem \ref{TSOFdh} in Subsection \ref{pt2}. In order to prove Theorem \ref{TSOF}, we first establish the relation between the kernels of $e^{it\Dk}$ and $e^{-it H_\K}$ then the proof follows from Theorem \ref{TSOFdh}. The following lemma will give the link between the kernels of $e^{it\Dk}$ and $e^{-it H_\K}$:
\begin{lem}\label{rel}
Let $K_{it}(x, y)$ and $L_{it}(x, y)$ be the  same as in \eqref{hsk} and \eqref{lsk} which are the  kernels of $e^{-it\Hk}$ and $e^{it\Dk}$ respectively. For $v>0$ we have \begin{equation}\label{kr}
K_{i\frac{\tan^{-1}v}{2}}(x, y)=(1+v^2)^{\frac{d+2\gk}{4}}e^{-\frac{iv}{2}|x|^2}L_{iv/2}(x\sqrt{1+v^2}, y).
\end{equation}
\end{lem}
\begin{proof}
The proof of this lemma can be abstracted from the proof of Proposition 8 in \cite{M3}. To make the article self contained we will only give the main steps. Let $v=\tan{2t}$ for $0<t<\frac{\pi}{4}$ then $t=\frac{1}{2}\tan^{-1}v$ and $\cot(2t)=\frac{1}{v}$ and $\sin(2t)=\frac{v}{\sqrt{1+v^2}}$ which implies
\begin{equation}\label{eq1}K_{\frac{i}{2}\tan^{-1}v}(x, y)=M_\K \Big(\frac{\sqrt{1+v^2}}{i v} \Big)^{\frac{d+2\gk}{2}} e^{\frac{i}{2v}(|x|^2+|y|^2)}E_\K\Big(\frac{x\sqrt{1+v^2}}{iv}, y\Big)\end{equation}
 and
\begin{equation}\label{eq2}L_{\frac{iv}{2}}(x\sqrt{1+v^2}, y)= M_\K \Big(\frac{1}{iv} \Big)^{\frac{d+2\gk}{2}} e^{\frac{i}{2v}\big((1+v^2)|x|^2+|y|^2\big)}E_\K\Big(\frac{x\sqrt{1+v^2}}{iv}, y\Big).\end{equation}
In view of \eqref{eq1} and \eqref{eq2} we get our required identity which completes the proof of the lemma.
\end{proof}
Now we are ready to prove Theorem \ref{TSOF}.
It can be verified using Lemma \ref{rel} and change of variable that for $p, q\geq 1$
\begin{equation}\label{id1} 
\left\|\sum_jn_j|e^{-itH_\K}f_j|^2\right\|^p_{L^p_tL^q_{\K, x}((0, \frac{\pi}{4})\times \R^d)} =\left \|\sum_jn_j|e^{-it\Dk}f_j|^2\right\|^p_{L^p_tL^q_{\K, x}((0, \infty)\times \R^d)}.
\end{equation}
For any orthonormal family of functions $(f_j)$ in $L^2_\K(\R^d)$ and $(n_j)\subset\C,$ let us define
\begin{equation*}
\varphi_{\{f_j\}}(t) := \left\|\sum_jn_j|e^{-itH_\K}f_j|^2\right\|_{L^q_{ \K}(\R^d)},
\end{equation*}
then the identity \eqref{id1} can be written as
 \begin{equation}\label{id2}\int^{\frac{\pi}{4}}_0|\varphi_{\{f_j\}}(t)|^p dt= \left\|\sum_jn_j|e^{-it\Dk}f_j|^2\right\|^p_{L^p_tL^q_{\K, x}((0, \infty)\times \R^d)}.\end{equation} 
Using the properties of kernel $K_{it}(x, y)$ given in \eqref{pk1} one can easily show that $\varphi_{\{f_j\}}$ satisfies the following:
\begin{equation}\label{rel2}
\varphi_{\{f_j\}}\Big(t+\frac{\pi}{2}\Big) = \varphi_{\{f_j\}}(t)~~ \text{and} ~~\varphi_{\{f_j\}}(-t) = \varphi_{\{\overline f_j\}}(t).
\end{equation}
Using \eqref{rel2} and \eqref{id2}, it can be easily proved that
\begin{equation}\label{id3}
\left(\int^\pi_{-\pi}|\varphi_{\{f_j\}}(t)|^p dt\right) = 4\int^{\frac{\pi}{4}}_{-\frac{\pi}{4}}|\varphi_{\{f_j\}}(t)|^p dt
=4 \int^\infty_{-\infty} \left\|\sum_jn_j|e^{-it\Dk}f_j|^2\right\|^p_{{L^q_\K}(\R^d)}dt.
\end{equation} 
Thus we establish
$$\left\|\sum_jn_j|e^{-itH_\K}f_j|^2\right\|^p_{L^p_tL^q_{\K, x}((-\pi, \pi)\times \R^d)}= 4\left \|\sum_jn_j|e^{-it\Dk}f_j|^2\right\|^p_{L^p_tL^q_{\K, x}(\R \times \R^d)}.$$
In view of the above identity, Theorem \ref{TSOF} follows from Theorem \ref{TSOFdh}.

\subsection{Proof of Theorem \ref{TSOFdh} and Theorem \ref{TSOFD}}\label{pt2}
We first prove Theorem \ref{TSOFdh} and then  Theorem \ref{TSOFD} follows by duality argument. In view of \eqref{id3} to prove \eqref{Ste3} it is enough to show that 
\begin{equation}\label{Ste33}
\int_{-\frac{\pi}{4}}^{\frac{\pi}{4}}\left(\int_{\R^d}\left|\sum_{j}n_j|\left(e^{-it\Hk}f_j\right)(x)|^2\right|^q dx\right)^{\frac{p}{q}}dt\leqslant{{C^{p}_{d,q}}\left(\sum_{j}|n_j|^{\frac{2q}{q+1}}\right)}^\frac{p(q+1)}{2q}.
\end{equation}
In order to prove inequality \eqref{Ste33} we use duality principle (Lemma \ref{Ld}) in view of which the inequality \eqref{Ste33} is equivalent to the estimate 
\begin{equation}\label{wts} 
 \|W_1T_{S}W_2\|_{\s^{\ap}(L^2_\K(\T\times\R^d))}\leq C\|W_1\|_{L^{\bt}_tL^{\ap}_{\K, x}(\T\times\R^d)} \|W_2\|_{L^{\bt}_tL^{\ap}_{\K, x}(\T\times\R^d)},
\end{equation}
where $\T=\big(-\frac{\pi}{4}, \frac{\pi}{4}\big),$ recall that $T_S=\ecs\ecs^*$ and $\ap,~~\bt$ satisfies 
\begin{equation}\label{conds}
\frac{2}{\bt}+\frac{d+2\gk}{\ap}=1,~~0\leq\frac{1}{\ap}<\frac{1}{d+2\gk+1}~~ \mbox{with}~~ \ap=\frac{2q}{q+1}.
\end{equation}
  Here $\ecs$ is considered as operator from $\ell^2(S)$ into $L^{p'}_tL^{q'}_{\K, x}(\T\times \R^d).$ It is easy to see that
$$T_Sf(t, x)= \int_{\R^d}\int_{-\frac{\pi}{4}}^{\frac{\pi}{4}}K(x, y, t-s)f(y, s) ds h^2_\K(y) dy,$$ where
\begin{equation}\label{kert}
K(x, y, t)= \sum_{\mu\in\Na_0^d}\phi_\mu^\K(x)\phi_\mu^\K(y)e^{-i(2|\mu|+d+2\gk)t}=K_{it}(x, y),
\end{equation}
where $K_{it}$ is defined as in \eqref{hsk}.  Inequality \eqref{Ste33} can be proved similar way as in Remark \ref{rem16} for $q=1$ and $p=\infty.$ This corresponds to inequality \eqref{wts} for $\ap=\infty.$ Therefore, with the help of interpolation theorem, it is enough to prove inequality \eqref{wts}, for $\ap$ such that $\frac{1}{d+2\gk+2}\leq \frac{1}{\ap}<\frac{1}{d+2\gk+1}.$ To prove this   
case we use Stein's complex interpolation theorem. In order to do that we define family of analytic operators $T_{z, \varepsilon}$ on the region $-1\leq\Re (z)\leq\frac{d+2\gk}{2}$ and $\varepsilon>0$ by 
$$T_{z, \varepsilon}f(t, x)= \int_{\R^d}\int_{-\frac{\pi}{4}}^{\frac{\pi}{4}}K_{z, \varepsilon}(x, y, t-s)f(y, s) ds h^2_\K(y) dy,$$
where $$K_{z, \varepsilon}(x, y, t)=t^z 1_{\varepsilon<|t|<\frac{\pi}{4}}K(x, y, t)$$ where $K(x, y, t)$ is defined as in \eqref{kert}. Since $K(x, y, t)=K_{it}(x, y)$ and in view of \eqref{hsk} we have that
\begin{equation}\label{kzee}
|K_{z, \varepsilon}(x, y, t)|= C |t|^{\Re (z)-\frac{(d+2\gk)}{2}},
\end{equation}
for $|t|<\pi/4$. An application of Hardy-Littlewood-Sobolev inequality (see \cite{BW}) along with inequality \eqref{kzee} yields
\begin{align*}
&\|W_1T_{z, \varepsilon}W_2\|^2_{\s^2(L^2_\K(\T\times\R^d))}\\
&=\int_{\T^2}\int_{\R^{2d}}|W_1(t, x)|^{2}|K_{z, \varepsilon}(x, y, t-s)|^2 |W_2(s, y)|^{2}h^2_\K(x) h^2_\K(y)  dx dy dt ds\\
&\leq C \int_{\T^2}\int_{\R^{2d}} \frac{|W_1(t, x)|^{2} |W_2(s, y)|^{2}}{|t-s|^{d+2\gk-2\Re (z)}}h^2_\K(x) h^2_\K(y)  dx dy dt ds\\
&\leq C \int_{\T}\int_{\T} \frac{\|W_1(t, \cdot)\|^{2}_{L^{2}_{\K, x}(\R^d)}{\|W_2(s, \cdot)\|^{2}_{L^{2}_{\K, x}(\R^d)}}}{|t-s|^{d+2\gk-2\Re (z)}} dt ds\\
&\leq C\left(\int_\T{\|W_1(t, \cdot)\|^{2\tilde{r}}_{L^{2}_{\K, x}(\R^d)}}dt\right)^\frac{1}{\tilde{r}} \left(\int_\T{\|W_2(t, \cdot)\|^{2\tilde{r}}_{L^{2}_{\K, x}(\R^d)}}dt\right)^\frac{1}{\tilde{r}},
\end{align*}
where $0\leq d+2\gk-2\Re (z)<1$ and $\frac{1}{\tilde r}+\frac{1}{2}(d+2\gk-2\Re (z))= 1.$
Let $r= 2\tilde{r}$ then we get that 
\begin{equation}\label{wtw}
\|W_1T_{z, \varepsilon}W_2\|_{\s^2(L^2_\K(\T\times\R^d))}\leq C {\|W_1\|_{L^r_tL^{2}_{\K, x}(\T\times\R^d)}} {\|W_2\|_{L^r_tL^{2}_{\K, x}(\T\times\R^d)}}
\end{equation}
provided $\frac{2}{r}+\frac{1}{2}(d+2\gk-2\Re (z))=1$ and $\Re (z)\in \Big(\frac{d+2\gk-1}{2}, \frac{d+2\gk}{2}\Big].$

Now we will prove that 
\begin{equation}\label{wtso}
\|W_1T_{z, \varepsilon}W_2\|_{\s^\infty(L^2_\K(\T\times\R^d))}\leq C(b) {\|W_1\|_{L^\infty L^\infty_{\K, x}(\T\times\R^d)}} {\|W_2\|_{L^\infty L^\infty_{\K, x}(\T\times\R^d)}},
\end{equation}
where $C(b)$ is a constant that grows exponentially in $b.$ Since $\s^\infty-$norm is the operator norm, to prove the inequality \eqref{wtso} it is enough to show that 
$$\|T_{-1+ib, \varepsilon}\|_{L^2_\K(\T\times\R^d)\to L^2_\K(\T\times\R^d)}\leq C(b),$$
equivalently 
$$\|T_{-1+ib, \ve}f\|_{L^2_\K(\T\times\R^d)}\leq C(b)\|f\|_{L^2_\K(\T\times\R^d)}.$$
Consider
\begin{align*} 
T_{-1+ib,\ve}f(t, x)&= \int_{\R^d}\int_{-\frac{\pi}{4}}^{\frac{\pi}{4}}K_{-1+ib, \ve}(x, y, t-s)f(y, s) ds h^2_\K(y) dy\\
&=\int_{\R^d}\int_{-\frac{\pi}{4}}^{\frac{\pi}{4}}(t-s)^{-1+ib} 1_{\varepsilon<|t|<\frac{\pi}{4}}(t-s)\\&\qquad\qquad\qquad\times\sum_{\mu\in\Na_0^d}\phi_\mu^\K(x)\phi_\mu^\K(y)e^{-(2|\mu|+d+2\gk)i(t-s)}f(s, y)ds h^2_\K(y) dy\\
&=\sum_{\mu\in\Na_0^d}\int_{-\frac{\pi}{4}}^{\frac{\pi}{4}}(t-s)^{-1+ib} 1_{\varepsilon<|t|<\frac{\pi}{4}}(t-s)\widehat{f_x}(s, \cdot)(\mu)\phi_\mu^\K(x)e^{-(2|\mu|+d+2\gk)i(t-s)}ds,
\end{align*}
where $\widehat{f_x}(s, \cdot)(\mu)= \int_{\R^d}f(s, y)\phi_\mu^\K(y)h^2_\K(y) dy,$
using Plancherel theorem for Dunkl-Hermite expansion, we get  
$$\|T_{-1+ib}f(t, \cdot)\|^2_{L^2_{\K, x}(\R^d)}= \sum_{\mu\in\Na_0^d}\left|\int_{-\frac{\pi}{4}}^{\frac{\pi}{4}}(t-s)^{-1+ib} 1_{\varepsilon<|t|<\frac{\pi}{4}}(t-s)G_\mu(s)ds\right|^2,$$
where $G_\mu(s)=\widehat{f_x}(s, \cdot)(\mu)e^{i(2|\mu|+d+2\gk)s}.$ 
Define 
$$F_b: G\mapsto \int_{-\frac{\pi}{4}}^{\frac{\pi}{4}}(t-s)^{-1+ib} 1_{\varepsilon<|t|<\frac{\pi}{4}}(t-s)G(s)ds.$$
With this notation we obtain that 
$$\|T_{-1+ib}f\|^2_{L^2_{\K}(\T\times\R^d)}= \sum_\mu\|F_bG_\mu(t)\|_{L^2(\T)}^2.$$ 
In view of the above and the identity $$\|f\|_{L^2_{\K}(\T\times\R^d)}^2=\sum_{\mu}\|G_\mu\|^2_{L^2(\T)},$$ it is enough to prove that $\|F_bG\|_{L^2(\T)}\leq C(b) \|G\|_{L^2(\T)}$ to get the required inequality \eqref{wtso}. With simple calculation it can be easily seen that 
\begin{equation}\label{eq.ht} F_bG(t)=\int_{\ve<|s|<\pi/4}\frac{h(t-s)}{s^{1-ib}} ds\end{equation} where
$$h(t)=\begin{cases} G(t), ~~&\mbox{when}~|t|<\pi/4\\ 0~~&\mbox{otherwise}.\end{cases}$$
Right hand side of \eqref{eq.ht} is just a Hilbert transform of $h$ up to $ib$, from \cite{VE}, this operator is bounded from $L^2(\R)$ into $L^2(\R)$ with the operator norm depends on $b$ exponentially. Thus we get $\|F_bG\|_{L^2(\T)}\leq C(b) \|G\|_{L^2(\T)}$.

In view of complex interpolation on the region $-1\leq\Re (z)\leq \la,$ where $\la\in \Big(\frac{d+2\gk-1}{2}, \frac{d+2\gk}{2}\Big],$ \eqref{wtw} and \eqref{wtso}, imply that 
$$\|W_1T_{0, \varepsilon}W_2\|_{\s^\ap(L^2_\K(\T\times\R^d))}\leq C\|W_1\|_{L^{\bt}_tL^{\ap}_{\K, x}(\T\times\R^d)} \|W_2\|_{L^{\bt}_tL^{\ap}_{\K, x}(\T\times\R^d)},$$
where $\ap$ and $\bt$ satisfies \eqref{conds}. By letting $\varepsilon\to 0,$ in the above, we get required inequality \eqref{wts}.   

\section{Applications of Strichartz estimates: Global well-posedness for the Dunkl Hermite-Hartree equation in Schatten spaces}
In this section we show the well-posedness results in the spirit of \cite{LS1, LS2} for a system of infinitely many equations (without a trace class assumption) in Schatten spaces for the  Dunkl Hermite-Hartree equation by applying our orthonormal Strichartz inequalities associated with Dunkl Hermite operator. Recall that $f\ast_\K g$ denotes the generalized (Dunkl) convolution of functions $f$ and $g$ (for the definition see in \cite{TX}).

\begin{theorem}
Let $$1\leq q<1+\frac{2}{d+2\gk-1}~~ \mbox{and}~ ~\frac{2}{p}+\frac{d+2\gk}{q} = d+2\gk$$
and $w: \R^d\rightarrow\C$ such that $f\mapsto w\ast_\K f$ is bounded operator from $L^p_\K(\R^d)\rightarrow L^\infty_\K(\R^d).$ Then for any $\gm_0\in\s^{\frac{2q}{q+1}},$ there exists a unique $\gm\in C^0_t([0, T], \s^{\frac{2q}{q+1}}(L^2_\K(\R^d))$ satisfying $\rho_\gm\in L^p_tL^q_{\K, x}([0, T] \times \R^d)$ and
\begin{align*} 
i\partial_t\gm&=[H_\K + w\ast_\K \rho_\gm, \gm],\\
&\gm|_{t=0}=\gm_0.
\end{align*}
\end{theorem} 
\begin{proof}
Let $R > 0$ such that $\|\gm_0\|_{\s^\frac{2q}{q+1}} = R<\infty.$ Let $T= T(R) (\leq 1$ to be chosen later). Consider the space 
$$X_T = \{(\gm, \rho) \in C^0_t([0, T], \s^{\frac{2q}{q+1}}(L^2_\K(\R^d))\times L^p_tL^q_{\K, x}([0, T] \times \R^d) : \|(\gm, \rho)\|_{X_T}\leq 4~ max\{1, C_{Stri}\}R\},$$ where the norm $\|(\gm, \rho)\|_{X_T}$ is defined by $$\|(\gm, \rho)\|_{X_T} := \|\gm\|_{ C^0_t([0, T], \s^{\frac{2q}{q+1}}(L^2_\K(\R^d))}+\|\rho\|_{L^p_tL^q_{\K, x}([0, T] \times \R^d)}.$$
Consider the map $$\Phi_1(\gm, \rho)(t) = e^{-itH_\K}\gm_0 e^{itH_\K} - i\int^t_0 e^{i(s-t)H_\K}  [w\ast_\K\rho_\gm(s), \gm(s)] e^{-i(s-t)H_\K} ds.$$ 
Using the above map, we define the contraction map $\Phi$ by $$\Phi(\gm, \rho) = (\Phi_1(\gm, \rho), \rho[\Phi_1(\gm, \rho)]),$$ where we have used the notation $\rho[\gm] = \rho_\gm.$ Now
\begin{align*}
 \|\Phi_1(\gm, \rho)\|_{ C^0_t([0, T], \s^{\frac{2q}{q+1}}(L^2_\K(\R^d))}& \leq R + 2 \int^T_0 \|w\ast_\K\rho(s)\|_{L^\infty_{\K, x} \|\gm(s)\|_{\s^\frac{2q}{q+1}}} ds\\
&\leq R + 2T^{1/p'} C_w \|\rho\|_{L^t_pL^q_{\K, x}}\|\gm\|_{C^0_t\s^\frac{2q}{q+1}}\\
& \leq R + 8T^{1/p'} C_w max (1, C^2_{Stri}) R^2
\end{align*}
and from Corollary \ref{ISI} we also have 
$$\|\rho[\Phi_1(\gm, \rho)]\|_{L^t_pL^q_{\K, x}} \leq C_{Stri} R + 8C_{Stri} T^{1/p'} C_w  max (1, C^2_{Stri}) R^2.$$
Choosing $T\leq 1$ small enough so that
$$8C_{Stri} T^{1/p'} C_w  max (1, C^2_{Stri}) R^2 \leq C_{Stri} R$$
and $\Phi$ maps $X$ to itself. Thus $\Phi$ is a contraction mapping and has a unique fixed point on $X$ which is a solution to the Hatree equation on $[0, T].$ 
\end{proof}


\section{Alternative proof of Theorem \ref{TSOF}}

\subsection{Proof of Theorem \ref{TSOF} and \ref{TSOFD} in the case $P=-\Dk$}
Throughout this section $\s^p$ stands for $\s^p(L^2_\K(\R^d))$.
As discussed before, the duality argument shows that Theorem \ref{TSOF} is equivalent to Theorem\ref{TSOFD} for the case $P=-\Dk$. By assuming $V\in L^\infty_c(\R\times\R^d)$ and that $\gamma$ is finite rank, we have to show that the operator  
$$ L^{p'}_tL^{q'}_{\K, x}(\R\times \R^d)\ni V\mapsto\int_\R e^{-it\D_\kappa}V(t, x)e^{it\D_\kappa}dt\in\s^{2q'},$$ is bounded. By using complex interpolation method \cite[Chap. 4]{BL},  it is enough to show the boundedness for the two points $(p', q')=(1, \infty)$ and $(p', q')=(1+\frac{d+2\gamma_\kappa}{2}, 1+\frac{d+2\gamma_\kappa}{2}).$ 

For the point  $(p', q')=(1, \infty)$: We will get with the well known argument that. 
\begin{align*}
\left\|\int_\R e^{-it\D_\kappa}V(t, x)e^{it\D_\kappa}dt\right\|_{\s^\infty}&\leq \int_\R \|e^{-it\D_\kappa}V(t, x)e^{it\D_\kappa}\|_{\s^\infty}dt\\ &=\int_\R\|V(t, \cdot)\|_{L^\infty_\K(\R^d)}dt=\|V\|_{L^1_tL^\infty_{\K, x}(\R \times \R^d)}.\end{align*} 

Now it is time to go for the case $p'= q'= 1+\frac{d+2\gamma_\kappa}{2}.$ Without any loss of generality, we may assume that $V\geq 0$ and  $V\in L^\infty_c(\R\times\R^d).$  Then we have $e^{-it\D_\kappa}V(t, x)e^{it\D_\kappa} \geq 0$ as an operator on $L^2_\K(\R^d),$ for all $t$ in $\R.$  It can be easily proved that $e^{-it\D_\kappa} x_j e^{it\D_\kappa}=x_j-2itT_j$ for each $j=1, 2, \ldots , d,$ using Dunkl transform (for the definition, see \cite{MR}) and its properties. So it can be written as $$e^{-it\D_\kappa} x e^{it\D_\kappa}=x-2it\nabla_\K,$$ where $\nabla_\K=(T_1, T_2, \ldots , T_d)$ and $x$ is identified with the multiplication operator.  By the functional calculus we deduce that 
\begin{equation}\label{alt1}
e^{-it\D_\kappa} f(x) e^{it\D_\kappa}=f(x+2tp),$$ where $p:=-i\nabla_\K.$ From this we can get $$e^{-it\D_\kappa} V(t, x) e^{it\D_\kappa}=V(t, x+2tp).
\end{equation}
Using the fact $V\geq 0,$ we can write the Schatten norm as 
\begin{align*}
&\left\|\int_\R e^{-it\D_\kappa}V(t, x)e^{it\D_\kappa}dt\right\|^{d+2\gk+2}_{\s^{d+2\gk+2}}\\
&=Tr\left(\int_\R e^{-it\D_\kappa}V(t, x)e^{it\D_\kappa}dt\right)^{d+2\gk+2}=Tr\left(\int_\R V(t, x+2tp)dt\right)^{d+2\gk+2}\\
&=Tr\left(\int_\R\cdot\cdot\cdot\int_\R V(t_1, x+2t_1p)\cdot\cdot\cdot V(t_{d+2\gk+2}, x+2t_{d+2\gk+2}p)\right)dt_1\cdot\cdot\cdot dt_{d+2\gk+2}.
\end{align*}

To exchange the trace and the integral, we need to prove that
\begin{equation}\label{TNE}
\int_\R\cdot\cdot\cdot\int_\R \|V(t_1, x+2t_1p)\cdot\cdot\cdot V(t_{d+2\gk+2}, x+2t_{d+2\gk+2}p)\|_{\s^1}dt_1\cdot\cdot\cdot dt_{d+2\gk+2}<\infty.
\end{equation}
(We are assuming $V\in L^\infty_c(\R\times\R^d))$ throughout). In order to estimate the Schatten norm in the integral, we make use of the following.
\begin{lem}\label{LE}
Let $\ap, \bt, \gm, \dl \in \R$. We have
\begin{equation}\label{l1}
\|f(\ap x +\bt p)g(\gm x+\dl p)\|_{\mathfrak{S}^r}\leq \frac{M_\K^{\frac{2}{r}}\|f\|_{L_\K^r(\R^d)}\|g\|_{L_\K^r(\R^d)}}{|\ap\dl -\bt \gm|^{(d+2\gm_\K)/r}},
\end{equation}
for all $r\geq 2$.
\end{lem}
 For the value $\alpha=\dl=1, \bt=\gm=0,$ and the multiplicity function $\K=0$, the above estimate is the well-known Kato-Seiler-Simon inequality (\cite{SS} and \cite[Thm. 4.1]{S}).    

We refer Subsection \ref{plsb} for the proof of this lemma. We now estimate the trace norm \eqref{TNE}. By using the fact $V\geq 0$ and H\"older's inequality in Schatten space, we can write
\begin{align*}
\|V(t_1, x+2t_1p)\cdots &V(t_{d+2\gk+2}, x+2t_{d+2\gk+2}p)\|_{\s^1}\\
&=\|V(t_1, x+2t_1p)\sqrt{V(t_2, x+2t_2p)}\sqrt{V(t_2, x+2t_2p)}\\ &\quad \times \cdot \cdot \cdot \times \sqrt{V(t_{d+2\gk+1}, x+2t_{d+2\gk+1}p)}V(t_{d+2\gk+2}, x+2t_{d+2\gk+2}p)\|_{\s^1}\\
&\leq\|V(t_1, x+2t_1p)\sqrt{V(t_2, x+2t_2p)}\|_{\s^{d+2\gk+1}}\\&\quad\|\sqrt{V(t_2, x+2t_2p)}\sqrt{V(t_3, x+2t_3p)}\|_{\s^{d+2\gk+1}}\\ &\quad\times\cdot\cdot\cdot \times\|\sqrt{V(t_{d+2\gk+1}, x+2t_{d+2\gk+1}p)}V(t_{d+2\gk+2}, x+2t_{d+2\gk+2}p)\|_{\s^{d+2\gk+1}}.
\end{align*}

Using \eqref{l1} and the fact $V\in L^\infty_c(\R\times\R^d),$ we get

\begin{align*}
\|&V(t_1, x+2t_1p)\cdot\cdot\cdot V(t_{d+2\gk+2}, x+2t_{d+2\gk+2}p)\|_{\s^1}\\ 
&\leq \frac{\|V(t_1, \cdot)\|_{L^{d+2\gk+1}_\K} \|V(t_2, \cdot)\|_{L^{({d+2\gk+1})/2}_\K}\cdot\cdot\cdot \|V(t_{d+2\gk+1}, \cdot)\|_{L^{({d+2\gk+1})/2}_\K}\|V(t_{d+2\gk+2}, \cdot)\|_{L^{d+2\gk+1}_\K}}{c^2_\K 2^{d+2\gk}|t_1-t_2|^{\frac{d+2\gk}{d+2\gk+1}} \cdot\cdot\cdot |t_{d+2\gk+1}-t_{d+2\gk+2}|^{\frac{d+2\gk}{d+2\gk+1}}}\\
&\leq C \frac{\prod^{d+2\gk+2}_{j=1}\mathbbm{1}(a\leq t_j\leq b)}{c^2_\K 2^{d+2\gk}|t_1-t_2|^{\frac{d+2\gk}{d+2\gk+1}} \cdot\cdot\cdot |t_{d+2\gk+1}-t_{d+2\gk+2}|^{\frac{d+2\gk}{d+2\gk+1}}},
\end{align*}
where $(a, b)$ is the support of $V$ in the time variable. Now we will use the multilinear Hardy-Littlewood-Sobolev inequality.

\begin{theorem}[Multilinear Hardy-Littlewood-Sobolev inequality  (See \cite{BW})]\label{HARDY}
Assume that $(\beta_{ij})_{1\leqslant i, j\leqslant N}$ and $(r_k)_{1\leqslant k\leqslant N}$ are real-numbers such that 
$$\beta_{ii}=0,~~ 0\leqslant\beta_{ij}=\beta_{ji}<1, r_k>1,~~ \sum^N_{k=1}\frac{1}{r_k}>1, ~~\sum^N_{i=1}\beta_{ik}=\frac{2(r_k-1)}{r_k}.$$
Then there exist a constant $C$ such that 
\begin{equation}\label{HINQ}
\left|\int_\R \cdots \int_\R\frac{f_1(t_1)\cdots f_N(t_N)}{\prod_{i<j}|t_i-t_j|^{\beta_{ij}}}dt_1\cdots dt_N\right|\leqslant C\prod^N_{k=1}\|f_k\|_{L^{r_k}(\R)},
\end{equation}
 for all $f_k\in L^{r_k}.$
 
\end{theorem}

The particular case of Multilinear Hardy-Littlewood-Sobolev inequality where all the $\bt_{ij}$ and the $r_k$ are identical, can be found in [\cite{CM}, Prop. 2.2]. Applying \eqref{HINQ} with taking $r_1=r_{d+2\gk+2}=\frac{2(d+2\gk+1)}{(d+2\gk+2)}$ and $r_2=\cdot \cdot \cdot = r_{d+2\gk+1}=d+2\gk+1,$ then we will get \eqref{TNE}. Hence now we can write
\begin{align*}     
Tr&\left(\int_\R e^{-it\D_\kappa}V(t, x)e^{it\D_\kappa}dt\right)^{d+2\gk+2}\\
&=\int_\R\cdot\cdot\cdot\int_\R Tr\left(V(t_1, x+2t_1p)\cdot\cdot\cdot V(t_{d+2\gk+2}, x+2t_{d+2\gk+2}p)\right)dt_1\cdot\cdot\cdot dt_{d+2\gk+2}.
\end{align*}

By the previous argument we will now obtain a more symmetric estimate of the trace in the integral. Simply using the cyclicity of the trace we get
\begin{align*}
&\left|Tr\left(V(t_1, x+2t_1p)\cdot \cdot \cdot V(t_{d+2\gk+2}, x+2t_{d+2\gk+2}p)\right)\right|\\ &\quad \quad \quad \quad
=\left|Tr\left(\sqrt {V(t_1, x+2t_1p)} \cdot \cdot \cdot V(t_{d+2\gk+2}, x+2t_{d+2\gk+2}p \sqrt{ V(t_1, x+2t_1p)}\right)\right|\\ &\quad \quad \quad \quad
\leq \|\sqrt {V(t_1, x+2t_1p)} \sqrt {V(t_2, x+2t_2p)}\|_{\s^{d+2\gk+2}}\\ &\quad \quad \quad \quad \quad \times \cdot \cdot \cdot \times 
\|\sqrt {V(t_{d+2\gk+1}, x+2t_{d+2\gk+1}p)} \sqrt {V(t_{d+2\gk+2}, x+2t_{d+2\gk+2}p)}\|_{\s^{d+2\gk+2}}\\ &\quad \quad \quad \quad \quad 
 \times \|\sqrt {V(t_{d+2\gk+2}, x+2t_{d+2\gk+2}p)} \sqrt {V(t_1, x+2t_1p)}\|_{\s^{d+2\gk+2}}.
\end{align*}
Again using of \eqref{l1} and multilinear Hardy-Littlewood-Sobolev inequality \eqref{HINQ} with $r_1 = \cdot \cdot \cdot = r_{d+2\gk+2} = 1+\frac{d+2\gk}{2}$ we will get our required inequality

$$Tr\left(\int_\R e^{-it\D_\kappa}V(t, x)e^{it\D_\kappa}dt\right)^{d+2\gk+2} \leq C \|V\|^{d+2\gk+2}_{L^{1+\frac{d+2\gk}{2}}_\K (\R\times \R^d)}.$$

\subsection{Proof of  Lemma \ref{LE}}\label{plsb}
 We prove the inequality \eqref{l1} for $r=2$ and $r=\infty$, then the general case follows from the complex interpolation. Since $e^{-it\D_\K}f(x)e^{it\D_\K}=f(x+2tp)$ with $p=-i\Nk$, we get that
$$e^{-it\D_\K}f( \ap x)e^{it\D_\K}=f( \ap x+2t\ap p).$$ If we substitute $t=\frac{\bt}{2\ap}$ in the above, we get that 
\begin{equation}f( \ap x+\bt p)=e^{-i\frac{\bt}{2\ap}\D_\K}f( \ap x)e^{i\frac{\bt}{2\ap}\D_\K},\end{equation} 
where $p=-i\N_\K.$ Similarly we have that 
\begin{equation}g( \gamma x+\delta p)=e^{-i\frac{\delta}{2\gamma}\D_\K}g( \gamma x)e^{i\frac{\delta}{2\gamma}\D_\K}.\end{equation} For $r=\infty,$ by noting that the Schatten $\infty$-norm is the operator norm on $L^2_\K(\R^d)$, the inequality \eqref{l1} can be seen easily, indeed
\begin{align*}
\|f(\ap x+\bt p)g(\gamma x+\delta p)\phi\|_{L^2_\K(\R^d)}&=\| e^{-i\frac{\bt}{2\ap}\D_\K}f( \ap x)e^{i\ap_0\D_\K}g( \gamma x)e^{i\frac{\delta}{2\gamma}\D_\K}\phi\|_{L^2_\K(\R^d)}\\
&\leq \|f\|_{L^\infty_\K(\R^d)}\|g\|_{L^\infty_\K(\R^d)}\|\phi\|_{L^2_\K(\R^d)},
\end{align*} where \begin{equation}\label{l2}\ap_0=\frac{\bt\gamma-\ap\delta}{2\ap\gamma}.\end{equation}
Here we used the facts $\|e^{it\D_\K}\phi\|_{L^2_\K(\R^d)}=\|\phi\|_{L^2_\K(\R^d)}$ and $\|f\phi\|_{L^2_\K(\R^d)}\leq\|f\|_{L^\infty_\K(\R^d)}\|\phi\|_{L^2_\K(\R^d)}$.  Hence $$\|f(\ap x+\bt p)g(\gamma x+\delta p)\|_{\s^\infty}\leq \|f\|_{L^\infty_\K(\R^d)}\|g\|_{L^\infty_\K(\R^d)}.$$
For $r=2,$ we know that $\|T\|_{\s^2}^2=Tr|T|^2$, we get
\begin{align*}
\|f(\ap x+\bt p)g(\gamma x+\delta p)\|_{\s^2}^2&=Tr[|f(\ap x+\bt p)|^2|g(\gamma x+\delta p)|^2]\\
&=Tr\Big[e^{-i\frac{\bt}{2\ap}\Dk}|f( \ap x)|^2e^{i\frac{\bt}{2\ap}\Dk}e^{-i\frac{\delta}{2\gamma}\Dk}|g( \gamma x)|^2e^{i\frac{\delta}{2\gamma}\Dk}\Big]\\
&=Tr\Big[|f( \ap x)|^2e^{i\ap_0\Dk}|g( \gamma x)|^2e^{-i\ap_0\Dk}\Big].
\end{align*}
We now show that the operator $|f( \ap x)|^2e^{i\ap_0\Dk}|g( \gamma x)|^2e^{-i\ap_0\Dk}$ is an integral operator with some kernel $k(x, y)$ using which we find the trace. Now
\begin{align*}
e^{it\D_\K}f(x) &=\int_{\R^d}L_{it}(x, y) f(y) h_\K^2(y) dy \\=&\frac{M_\K}{ (2it)^{\gk+d/2}}\int_{\R^d}e^\frac{{i(|x|^2+|y|^2})}{4t}E_\K\left(\frac{x}{2it}, y\right)f(y)h^2_\K(y)dy.
\end{align*}
So, for a suitable $\phi$ we have
$$|g(\gamma\xi)|^2e^{i\ap_0\D_\K}\phi(\xi)=\frac{M_\K}{ (2i\ap_0)^{\gk+d/2}}\int_{\R^d}e^\frac{{i(|\xi|^2+|y|^2)}}{4\ap_0}E_\K\Big(\frac{\xi}{2i\ap_0}, y\Big)|g(\gamma\xi)|^2\phi(y)h^2_\K(y)dy,$$
which implies
\begin{align*}
|f(\ap x)|^2e^{-i\ap_0\D_\K}|g(\gamma\cdot)|^2e^{i\ap_0\D_\K}\phi(x)=&\frac{M_\K^2}{ (2\ap_0)^{d+2\gk}}\int_{\R^d}\int_{\R^d}\left(e^\frac{-i(|x|^2+|\xi|^2)}{4\ap_0}e^\frac{i(|\xi|^2+|y|^2)}{4\ap_0}E_\K\Big(\frac{-x}{2i\ap_0}, \xi\Big) \right.\\ &\quad \times \left. E_\K\Big(\frac{\xi}{2i\ap_0}, y\Big)|f(\ap x)|^2|g(\gamma\xi)|^2\phi(y)  h^2_\K(y)h^2_\K(\xi)\right) dy d\xi\\
=&\int_{\R^d}k(x,y)\phi(y)h^2_\K(y)dy,
\end{align*}
where $$k(x, y)=\frac{M_\K^2}{(2\ap_0)^{d+2\gk}}\int_{\R^d}e^\frac{-i(|x|^2-|y|^2)}{4\ap_0} E_\K\Big(\frac{-x}{2i\ap_0}, \xi\Big)E_\K\Big(\frac{\xi}{2i\ap_0}, y\Big)|f(\ap x)|^2|g(\gamma\xi)|^2h^2_\K(\xi)d\xi.$$
Therefore,
$$k(x, x)=\frac{M_\K^2}{(2\ap_0)^{d+2\gk}}\int_{\R^d}E_\K\Big(\frac{-x}{2i\ap_0}, \xi\Big)E_\K\Big(\frac{\xi}{2i\ap_0}, x \Big)|f(\ap x)|^2|g(\gamma\xi)|^2h^2_\K(\xi)d\xi.$$

Now we can find the required trace as
\begin{align*}
Tr\big[|f( \ap x)|^2e^{i\ap_0\Dk}|g( \gamma x)|^2e^{-i\ap_0\Dk}\big]=&\int_{\R^d}k(x,x)h^2_\K(x)dx\\
\leq &\frac{M_\K^2}{(2|\ap_0|)^{d+2\gk}}\int_{\R^d}\int_{\R^d}|f(\ap x)|^2|g(\gamma\xi)|^2 h^2_\K(\xi) h^2_\K(x) d\xi dx\\
=&\frac{M_\K^2}{(2|\ap_0 \ap\gm|)^{d+2\gk}}\|f\|^2_{L^2_\K(\R^d)}\|g\|^2_{L^2_\K(\R^d)}.
\end{align*}
The first inequality in the above follows in view of the expression for $k(x, x)$ and the fact $|E_\K(ix, y)|\leq1$.
In view of \eqref{l2}, we have that $2|\ap_0 \ap \gm|=|\bt\gm-\ap \dl|$ and hence
$$\|f(\ap x+\bt p)g(\gamma x+\delta p)\|_{\s^2}^2\leq \frac{M_\K^2}{|\bt\gamma-\ap\dl|^{d+2\gamma_\K}}\|f\|^2_{L^2_\K}\|g\|^2_{L^2_\K}.$$
The inequality in $\s^r$ now follows from the complex interpolation.

\section{Final Remarks}
In this paper, using a relation between the kernels  of Schr\"odinger propagators  $e^{-it H_\kappa}$ and $e^{it\Delta_\kappa}$, we proved Strichartz inequality   associated with Dunkl Hermite operator $H_\kappa$ is for   a system of orthonormal functions with initial data.  However,   we are going to consider the following problems as the  future work:
\begin{enumerate}
\item We will try to prove that the range of the exponents given in Theorem \ref{TSOF} and Theorem \ref{TSOFdh} are necessary to satisfy the estimates \eqref{Ste2} and \eqref{Ste3}.

\item We will try to prove the Schatten exponent $\alpha=\frac{2q}{q+1}$ appear in the right hand side of the estimates \eqref{Ste2} and \eqref{Ste3} is optimal.
\end{enumerate} 

%

\section*{Acknowledgments}
Both the authors are grateful to the Department of Mathematics, BITS Pilani K K Birla Goa Campus for the facilities utilised. The authors thank Prof. S. Thangavelu, IISc Bangalore for the fruitful discussions.


\end{document}